\documentclass{article}
\usepackage{amssymb,amsmath,amsthm,graphicx,cite}
\usepackage{enumerate}
\newcommand{\Z}{\mathbb{Z}}

\def\utr{\, \underline{\triangleright}\, }
\def\otr{\, \overline{\triangleright}\, }

\newtheorem{theorem}{Theorem}

\newtheorem{corollary}[theorem]{Corollary}

\theoremstyle{definition}
\newtheorem{example}{Example}
\newtheorem{definition}{Definition}
\newtheorem{remark}{Remark}

\title{\Large \textbf{Biquandle Bracket Quivers}}

\date{}

\author{Pia Cosma Falkenburg\footnote{Email: ich@picofa.de}
\and
Sam Nelson\footnote{Email: Sam.Nelson@cmc.edu. Partially supported by Simons Foundation collaboration grant 702597.}}

\begin{document}
\maketitle

\begin{abstract}
Biquandle brackets define invariants of classical and virtual knots and links
using skein invariants of biquandle-colored knots and links. Biquandle coloring
quivers categorify the biquandle counting invariant in the sense of defining
quiver-valued enhancements which decategorify to the counting invariant. In this
paper we unite the two ideas to define biquandle bracket quivers, providing new
categorifications of biquandle brackets. In particular, our construction
provides an infinite family of categorifications of the Jones polynomial and
other classical skein invariants.
\end{abstract}

\parbox{6in} {\textsc{Keywords:} Biquandles, Biquandle Brackets, Quiver Enhancements, Enhancements of Counting Invariants

\smallskip

}

\section{\large\textbf{Introduction}}\label{I}

In \cite{NOR}, skein invariants of biquandle-colored knots and links known
as \textit{biquandle brackets} were introduced. A biquandle bracket associates
an element of a commutative ring with identity $R$ to a biquandle-colored knot 
or link diagram via skein relations which vary depending on the biquandle 
colors at the crossing, in a way which is unchanged by biquandle-colored
Reidemeister moves. The multiset of such ring elements over the complete set
of biquandle colorings (i.e, over the biquandle homset 
$\mathrm{Hom}(\mathcal{B}(L),X)$ where $\mathcal{B}(L)$ is the fundamental 
biquandle of the link $L$ and $X$ is the coloring biquandle) forms an enhanced
invariant of knots and links that determines the biquandle counting invariant
and includes the classical skein invariants and the biquandle 2-cocycle
invariant as special cases.

In \cite{CN}, the quandle homset invariant was categorified (in the sense of
defining a category-valued invariant from which the original invariant
can be recovered) to directed graphs we call \textit{quandle 
coloring quivers} associated to subsets of the endomorphism ring of the
coloring quandle $X$. Convenient polynomial enhancements of the quandle
counting invariant were obtained as decategorifications of these quivers.
In \cite{CCN} these quiver-valued invariants were extended to the cases of 
pseudoknots, singular knots and virtual knots via psyquandles and biquandles;
in particular, biquandle coloring quivers were established for classical and 
virtual knots and links.

In this paper we categorify biquandle brackets via biquandle counting quivers 
to obtain \textit{biquandle bracket quivers}, quiver-valued invariants of 
classical and virtual knots and links weighted with weights associated to
a subset of the endomorphism ring of the coloring biquandle together with a
biquandle bracket structure with coefficients in a commutative ring with
identity $R$. The original biquandle 
bracket polynomial invariants as well as new enhanced invariants are recovered 
as decategorifications. We note that this construction yields an infinite 
family of alternative categorifications of the classical skein invariants and 
the biquandle 2-cocycle invariants, one for each element of the power set of
each biquandle in the biquandle bracket construction.

The paper is organized as follows. In Section \ref{BBaQ} we review the basics
of biquandles, biquandle brackets and biquandle counting quivers. In Section 
\ref{BBQ} we introduce the new invariants and illustrate their computation
via examples. In Section \ref{E} we provide additional computational examples, 
and we conclude in Section \ref{Q} with some questions for future research.

\section{\large\textbf{Biquandles, Brackets and Quivers}}\label{BBaQ}

In this section we briefly review biquandles, biquandle brackets and quivers; 
see \cite{CN}, \cite{EN} and \cite{NOR} for more details.

\begin{definition}
A \textit{biquandle} is a set $X$ with two binary operations 
$\utr, \otr : X \rightarrow X$ satisfying $\forall x,y,z \in X$:
\begin{enumerate}[(i)]
\item $x \utr x = x \otr x$,
\item The maps $ \alpha_y, \beta_y: X \rightarrow X$ and 
$S: X \times X \rightarrow X \times X$ defined by $\alpha_y(x)=x \otr y$, 
$\beta_y(x)=x \utr y$ and $S(x,y)=(y \otr x, x \utr y)$ are invertible,
\item The \textit{exchange laws} are satisfied:
\begin{align*}
\left(x \utr y\right)\utr \left(z \utr y\right) &= \left(x \utr z\right)\utr \left(y \otr z\right) \\
\left(x \utr y\right)\otr \left(z \utr y\right) &= \left(x \otr z\right)\utr \left(y \otr z\right) \\
\left(x \otr y\right)\otr \left(z \otr y\right) &= \left(x \otr z\right)\otr \left(y \utr z\right)
\end{align*}
\end{enumerate}
If $x \otr y = x$ for all $x,y \in X$, we say $X$ is a \textit{quandle}.
\end{definition}

\begin{example}
Let $X=\Z_n$, $x \otr y = x$ for all $x,y \in X$ and $x \utr y \equiv 2y-x$ 
(mod $n$). Then $X$ is called a \textit{dihedral quandle}, also sometimes 
called a \textit{cyclic quandle} or a \textit{Takasaki kei}.
\end{example}

\begin{example}
Let $G$ be a group, $x \otr y = x$ for all $x,y \in X$ and let $\utr$ be 
$n$-fold conjugation: $x \utr y = y^nxy^{-n}$. Then G is a quandle.
\end{example}

\begin{example}
For a non-quandle family of examples, let $R$ be a commutative ring with
identity and let $t,s\in R^{\times}$ be units in $R$. Then $R$ has a biquandle
structure known as an \textit{Alexander biquandle} defined by
\[x\utr y = tx+(s-t)y\quad \mathrm{and}\quad x\otr y= sx.\]
\end{example}

\begin{example}
Let $L$ be an oriented knot or link. The fundamental biquandle of $L$, denoted $\mathcal{B}(L)$, is the set of equivalence classes of biquandle words in a set of generators corresponding to the semiarcs in a diagram of $L$ under the equivalence relation generated by the crossing relations of $L$ and the biquandle axioms. The Hopf link $L2a1$
\[\includegraphics{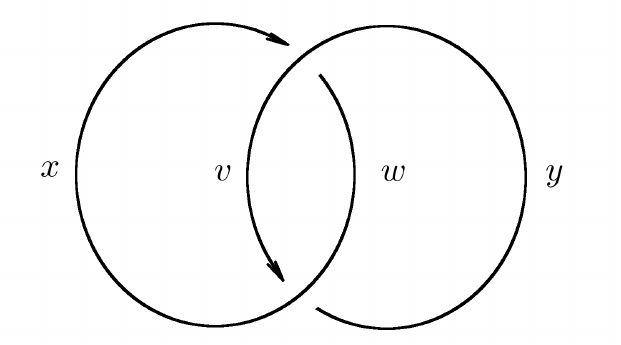}\]
has the fundamental biquandle presentation 
\[\mathcal{B}(L2a1)= \langle x,y,v,w \,|\, x \otr v = w, v \utr x = y, v \otr x = y, x\utr y = w \rangle\]
\end{example}

\begin{remark}
The fundamental biquandle of a classical knot or link interpreted as an 
Alexander biquandle determines the Alexander polynomials of the knot or link,
for virtual knots and links determines the Sawollek polynomial, also 
called the generalized Alexander polynomial. See \cite{KR} for more details.
\end{remark}

\begin{definition}
Let $X,Y$ be biquandles with binary operations $\otr_X,\utr_X,\otr_Y,\utr_Y$ respectively. A \textit{biquandle homomorphism} is a map $f: X \rightarrow Y$ such that for any $a,b \in X$ we have
\begin{align*}
f(a \otr_X \, b)= f(a) \otr_Y f(b)\\
f(a \utr_X \, b)= f(a) \utr_Y f(b)
\end{align*}
\end{definition}

\begin{definition}
Let $X$ be a finite biquandle, called the \textit{coloring biquandle} and 
$L$ be an oriented knot or link diagram. Then a \textit{biquandle coloring} of 
$L$ is an assignment of elements of $X$ to the semiarcs in $L$ such that the 
following crossing relations are satisfied at every crossing: 
\[\includegraphics{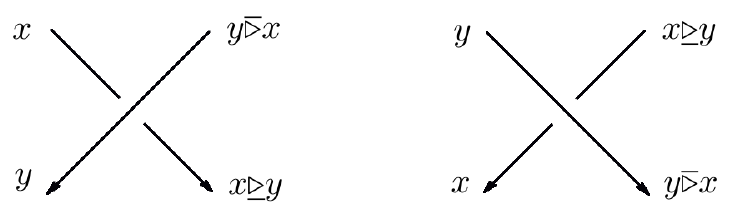}\]
\end{definition}

\begin{remark}
The set of biquandle colorings of $L$ can be identified with the set of 
biquandle homomorphisms from the fundamental biquandle of $L$ to $X$, denoted 
$\mathrm{Hom}(\mathcal{B}(L), X)$. 
If $L$ is tame and $n$ denotes the number of semiarcs in $L$ then $\mathcal{B}$
is finitely generated with $2n$ generators; hence 
$|\mathrm{Hom}(\mathcal{B}(L), X)| \leq |X|^{2n}$. The \textit{biquandle 
counting invariant} is the cardinality of the coloring space, denoted 
$\Phi^{\Z}_X (L) = |\mathrm{Hom}(\mathcal{B}(L), X)|$.
\end{remark}

\begin{definition}
Let $X$ be a finite biquandle, $R$ a commutative ring with identity and let $R^X$ denote the units of $R$. A \textit{biquandle bracket} over $R$ is a pair of maps $A,B : X\times X \rightarrow R^X$ assigning $A_{x,y}, B_{x,y} \in R^X$ to the pairs of elements $\left(x, y\right) \in X \times X$ such that the following conditions are satisfied:
\begin{enumerate}[(i)]
\item $\forall x \in X$, the elements $-A^2_{x,x}B^{-1}_{x,x}$ are equal with their common value denoted as $w \in R^X$, i. e. $w$ is the value of a positive crossing
\item $\forall x,y \in X$, the elements $-A^{-1}_{x,y}B_{x,y}-A_{x,y}B^{-1}_{x,y}$ are equal with their common value denoted as $\delta \in R$, i. e. $\delta$ is the value of a simple closed curve
\item $\forall x,y,z \in X$ we have
\begin{align*}
A_{x,y}A_{y,z}A_{x\utr y,z \otr y}=&A_{x,z}A_{y \otr x,z \otr x}A_{x\utr z,y\utr z}\\
A_{x,y}B_{y,z}B_{x\utr y,z \otr y}=&B_{x,z}B_{y \otr x,z \otr x}A_{x\utr z,y\utr z}\\
B_{x,y}A_{y,z}B_{x\utr y,z \otr y}=&B_{x,z}A_{y \otr x,z \otr x}B_{x\utr z,y\utr z}\\
A_{x,y}A_{y,z}B_{x\utr y,z \otr y}=&A_{x,z}B_{y \otr x,z \otr x}A_{x\utr z,y\utr z} + A_{x,z}A_{y \otr x,z \otr x}B_{x\utr z,y\utr z} \\
&+ \delta A_{x,z}B_{y \otr x,z \otr x}B_{x\utr z,y\utr z} + B_{x,z}B_{y \otr x,z \otr x}B_{x\utr z,y\utr z}\\
B_{x,y}A_{x,y}A_{x\utr y,z \otr y}+A_{x,y}B_{y,z}A_{x\utr y,z \otr y}\; \; \; \;&\\
+\delta B_{x,y}B_{y,z}A_{x\utr y,z \otr y}+B_{x,y}B_{y,z}B_{x\utr y,z \otr y} =&B_{x,z}A_{y \otr x,z \otr x}A_{x\utr z,y\utr z}
\end{align*}
\end{enumerate}
These biquandle bracket axioms are chosen in such a way that the state-sum 
expansion of a $X$-colored oriented knot or link diagram L using the skein 
relations 
\[\includegraphics{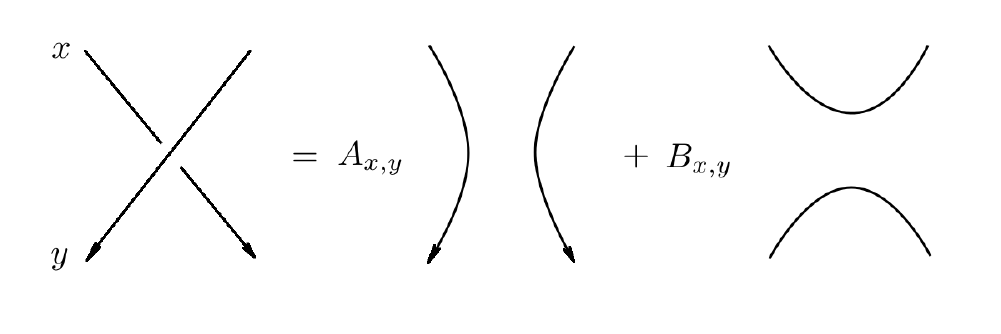}\]
\[\includegraphics{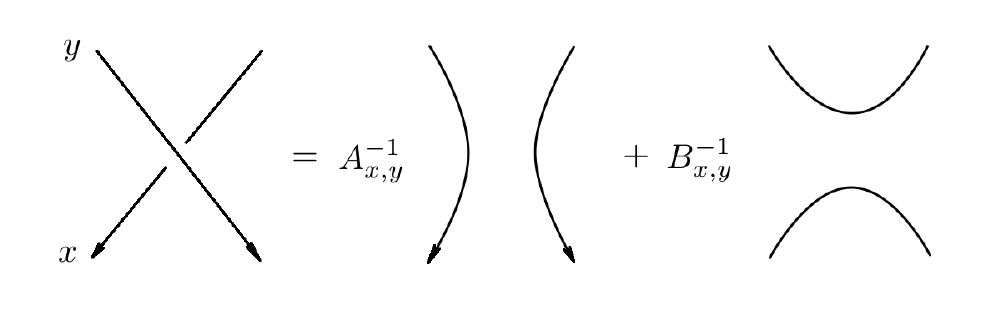}\]
is invariant under $X$-colored Reidemeister moves with $w$ as the value 
of a positive crossing and $\delta$ as the value of a simple closed curve. 
\[\includegraphics{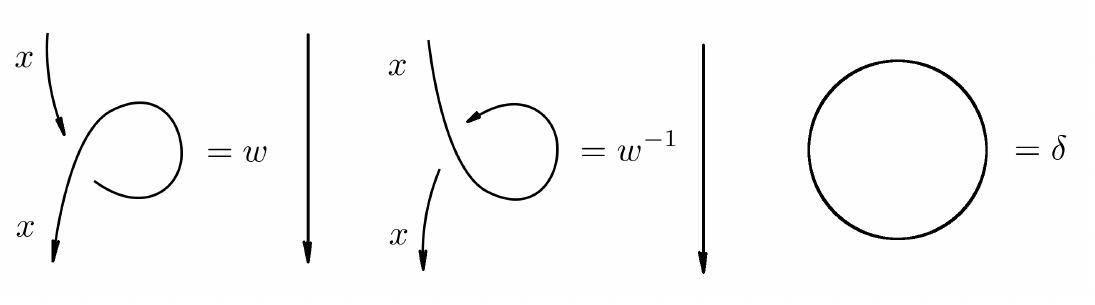}\]
\end{definition}

We specify biquandle brackets over $X=\{x_1,\dots, x_n\}$ with a block matrix
$[A|B]$ with blocks $A$ and $B$ whose $(j,k)$ entries are the skein
coefficients $A_{j,k}$ and $B_{j,k}$.

Examples of biquandle brackets include:
\begin{itemize}
\item Classical skein invariants including the Conway, Jones, HOMFLYPT and 
Kauffman polynomials with trivial biquandle,
\item Biquandle 2-cocycle invariants (see \cite{CES,EN}) are biquandle 
brackets such that the matrices of coefficients $A$ and $B$ are equal, and
\item Other potentially new invariants.
\end{itemize}

\begin{remark}\label{cohom}
Biquandle brackets representing biquandle 2-coboundaries can be used to modify
other biquandle brackets over the same biquandle and ring by entrywise 
multiplication to obtain new biquandle brackets, with the resulting invariant
equivalent to the original. Curiously, this operation extends the notion of 
biquandle brackets being cohomologous to other brackets beyond the set of
biquandle cocycles. In particular, in \cite{HVW} it is shown that many examples
in other works are cohomologous to brackets representing the Jones polynomial.
\end{remark}

\begin{definition}
Let $X$ be a finite biquandle and $L$ an oriented knot or link. For any set 
of biquandle endomorphisms $S \subset\mathrm{Hom}(X, X)$ the associated 
\textit{biquandle coloring quiver}, denoted $\mathcal{Q}^S_X (L)$, is the 
directed graph whose vertices are all elements 
$f \in\mathrm{Hom}(\mathcal{B}(L), X)$. An edge is directed from $f$ to $g$ 
if $g=\phi f$ for an element $\phi \in S$.
\end{definition}

An important special case is $S=\mathrm{Hom}(X, X)$, which is called the \textit{full biquandle coloring quiver} of $L$ with respect to $X$, denoted $\mathcal{Q}_X (L)$.
\begin{example}
The Hopf link $L2a1$ has four colorings by the biquandle $X$ given by the operation tables
\[\begin{array}{r|rrrr}
\utr & 1 & 2 \\ \hline
1 & 1 & 1\\
2 & 2 & 2  
\end{array}
\quad
\begin{array}{r|rrrr}
\otr & 1 & 2 & \\ \hline
1 & 1 & 1  \\
2 & 2 & 2
\end{array}\]
Choosing $S=\{\phi_1,\phi_2\}$ where $\phi_1(1)=\phi_1(2)=1$ and $\phi_2(1)=\phi_2(2)=2$ yields the following biquandle coloring quiver: 
\[\includegraphics{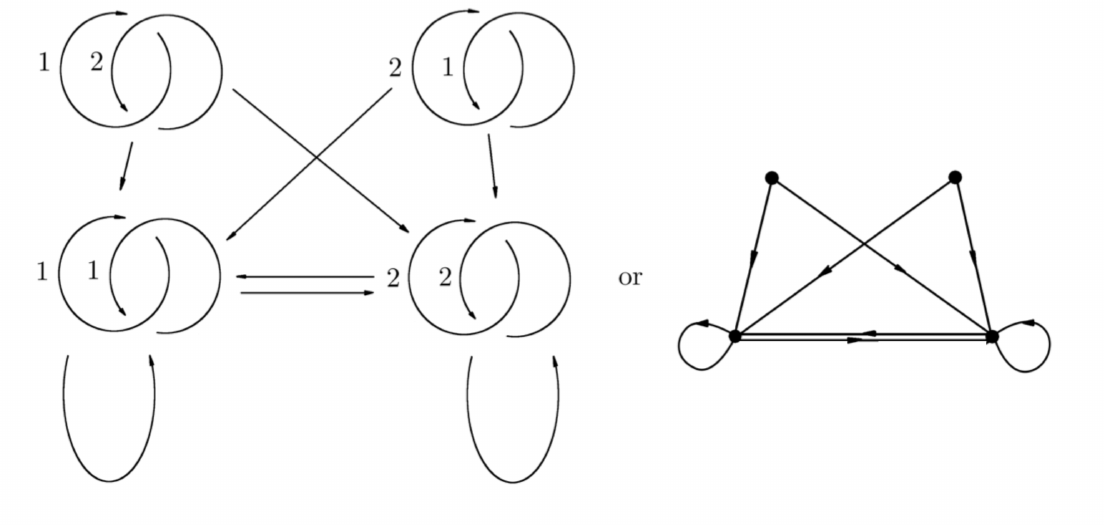}\]
\end{example}

Biquandle coloring quivers provide categorifications of the biquandle counting
invariant in the sense that they are quivers and hence small categories, with
vertices as objects and directed paths as morphisms (with empty paths based at
vertices as identity elements). Any further invariants of quivers then 
automatically provide invariants of knots and links via this construction.

\section{\large\textbf{Biquandle Bracket Quivers}}\label{BBQ}

In this section we introduce our main definition, the biquandle bracket quiver,
and illustrate its computation with an example.

\begin{definition}
Let $X$ be a finite biquandle, $R$ a commutative ring with identity,
$\beta=\{A,B:X\times X\to R^{\times}\}$ a biquandle bracket and 
$S\subset\mathrm{Hom}(X,X)$ a set of biquandle endomorphisms. 
For any oriented
link $L$, we define the \textit{biquandle bracket quiver} of $L$ with respect
to $X,R,\beta$ and $\phi$, denoted $\mathcal{BBQ}_{X,\beta}^S(L)$ to be the 
biquandle coloring quiver of $L$ with vertices $v$ weighted with the biquandle 
bracket values $\beta(v)$ of the $X$-colorings of $L$ they represent.
\end{definition}

Since both the biquandle coloring quiver and the biquandle bracket multiset
are unchanged by Reidmeister moves, we have our main result:

\begin{theorem}
The isomorphism class of $\mathcal{BBQ}_{X,\beta}^S(L)$ as a weighted graph
is an invariant of oriented knots and links.
\end{theorem}

\begin{example}\label{ex1}
Let $X$ be the biquandle given by the operation tables
\[
\begin{array}{r|rrrr}
\utr & 1 & 2 & 3 & 4 \\ \hline
1 & 1 & 1 & 2 & 2 \\
2 & 2 & 2 & 1 & 1 \\
3 & 3 & 3 & 4 & 4 \\
4 & 4 & 4 & 3 & 3
\end{array}
\quad
\begin{array}{r|rrrr}
\otr & 1 & 2 & 3 & 4 \\ \hline
1 & 1 & 1 & 1 & 1 \\
2 & 2 & 2 & 2 & 2 \\
3 & 4 & 4 & 4 & 4 \\
4 & 3 & 3 & 3 & 3,
\end{array}
\]
let $\beta$ be the biquandle bracket over $X$ with coefficients in 
$\mathbb{Z}_3$ given by 
\[\left[\begin{array}{rrrr|rrrr}
1 & 2 & 2 & 2 & 2 & 1 & 1 & 1 \\
1 & 1 & 1 & 2 & 2 & 2 & 2 & 1 \\
2 & 1 & 1 & 1 & 1 & 2 & 2 & 2 \\
2 & 2 & 2 & 1 & 1 & 1 & 1 & 2
\end{array}
\right],\]
which has $\delta=-(1)(2)-2(1)=-4=-1=2$ and $w=-1^2(2)=-2=1$.
let $\phi:X\to X$ be the biquandle endomorphism defined by $\phi(1)=2$, 
$\phi(2)=1$, $\phi(3)=3$ and $\phi(4)=4$ and consider the Hopf link $L=L2a1$. 
$L$ has eight $X$-colorings determining a biquandle coloring quiver as shown
\[\includegraphics{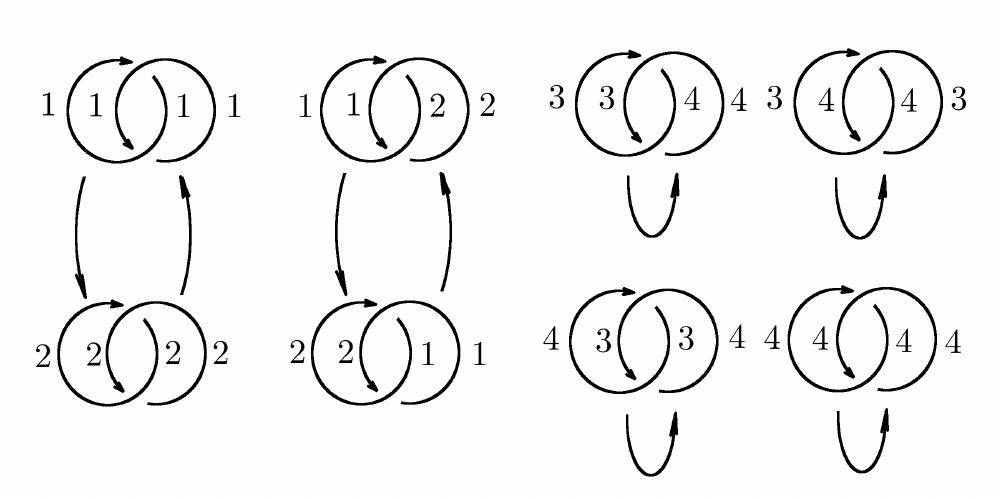}\]
Taking one of these colorings and expanding to the set of Kauffman states, 
we find the biquandle bracket value for this coloring
\[\includegraphics{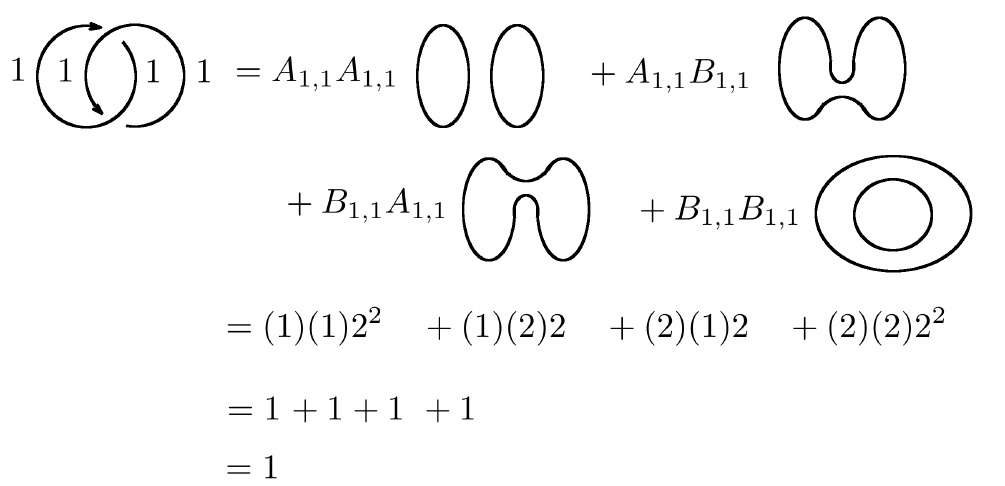}\]
Repeating for each coloring, we obtain the biquandle bracket quiver
$\mathcal{BBQ}_{X,\beta}^S(L)$
\[\includegraphics{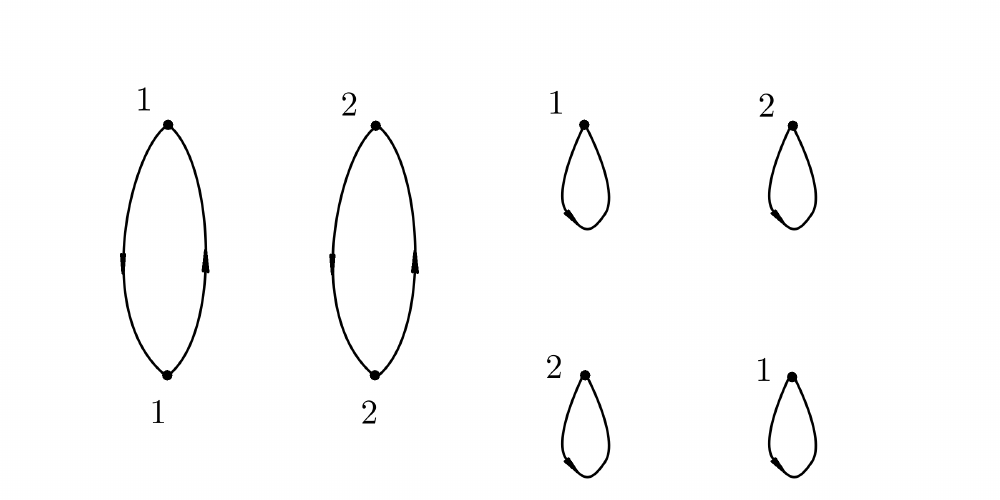}\]
\end{example}

As with the biquandle coloring quiver, $\mathcal{BBQ}_{X,\beta}^S(L)$ is a small 
category with vertices 
as objects and directed paths as morphisms (with identity morphisms given by 
``empty'' paths which stay at a vertex); hence, it is a categorification of
the biquandle bracket polynomial invariant. Invariants of these graphs then
give us new invariants of knots as decategorifications:

\begin{definition}
Let $\mathcal{BBQ}_{X,\beta}^S(L)$ be the biquandle bracket quiver of an
oriented link $L$ with vertex set $V$ and edge set $E$. We define the
\textit{in-degree biquandle bracket polynomial} to be
\[\Phi_{X,\beta}^{S,\mathrm{deg}^+}(L)=\sum_{v\in V} u^{\beta(v)}v^{\mathrm{deg}^+(v)}\]
and we define the \textit{biquandle bracket 2-variable polynomial} to be
\[\Phi_{X,\beta}^{S,2}(L)=\sum_{e\in E} s^{\beta(s(e))}t^{\beta(t(e))}\]
where $s(e)$ and $t(e)$ are the source and target vertices of each edge 
$e\in E$.
\end{definition}

\begin{corollary}
The polynomials $\Phi_{X,\beta}^{S,\mathrm{deg}^+}(L)$ and 
$\Phi_{X,\beta}^{S,\mathrm{deg}^+}(L)$ are invariants of oriented classical
knots and links.
\end{corollary}

\begin{remark}
If $S=\{\phi\}$ is a singleton, we will often write $\phi$ instead of 
$\{\phi\}$ in the superscript for simplicity. 
\end{remark}

\begin{example}
The Hopf link $L2a1$ has in-degree biquandle bracket and 2-variable biquandle 
bracket polynomials
\[\Phi_{X,\beta}^{\phi,\mathrm{deg}^+}(L2a1)=4u^2v+4uv\]
and
\[\Phi_{X,\beta}^{\phi,2}(L2a1)=4s^2t^2+4st\]
with respect to the biquandle $X$, biquandle bracket $\beta$ and biquandle
endomorphism $\phi$ in Example \ref{ex1}.
\end{example}

Everything in our construction applies without modification to virtual knots 
and links by the usual convention of ignoring the virtual crossings when 
determining biquandle colorings; hence we have the following:

\begin{corollary}
Biquandle bracket quivers $\mathcal{BBQ}_{X,\beta}^{S}(L)$ and the polynomials 
$\Phi_{X,\beta}^{\phi,2}(L)$ and $\Phi_{X,\beta}^{\phi,\mathrm{deg}^+}(L)$
are invariants of oriented virtual knots and links.
\end{corollary}

\section{\large\textbf{Examples}}\label{E}

In this section we collect a few examples of the new invariants.

\begin{example}
Let $X$ and $\beta$ be the biquandle and biquandle bracket on $R=\mathbb{Z}_6$ 
defined by the operation tables and matrix
\[
\begin{array}{r|rrrr}
\utr & 1 & 2 & 3 & 4 \\ \hline
1 & 2 & 2 & 2 & 2 \\
2 & 1 & 1 & 1 & 1 \\
3 & 4 & 4 & 4 & 4 \\
4 & 3 & 3 & 3 & 3
\end{array}\quad
\begin{array}{r|rrrr}
\otr & 1 & 2 & 3 & 4 \\ \hline
1 & 2 & 2 & 1 & 1 \\
2 & 1 & 1 & 2 & 2 \\
3 & 3 & 3 & 4 & 4 \\
4 & 4 & 4 & 3 & 3
\end{array}
\quad
\left[\begin{array}{rrrr|rrrr}
1 & 5 & 1 & 1 & 5 & 5 & 1 & 5 \\
5 & 1 & 1 & 5 & 1 & 5 & 5 & 1 \\
1 & 5 & 1 & 5 & 5 & 1 & 5 & 1 \\
1 & 5 & 5 & 1 & 5 & 1 & 1 & 5
\end{array}\right]
\] 
and let $\phi_1,\phi_2:X\to X$ be the endomorphism defined by 
$\phi(1)=2, \phi(2)=1, \phi(3)=4$ and $\phi(4)=3$. We computed via 
\texttt{python} the $\Phi_{X,\beta}^{\phi,2}$-values 
defined in Section \ref{BBQ} for a choice of orientation for each of the prime 
links with up to seven crossings in the table at \cite{KA}; the results are 
in the table.
\[
\begin{array}{r|l}
\Phi_{X,\beta}^{\phi_1,2}(L) & L \\ \hline
6s^4t^4+2s^2t^2 & L2a1, L6a2, L6a3, L7a5, L7a6,  \\
8s^4t^4+8s^2t^2 & L4a1 \\
8s^4t^4+2s^4t^2+2s^2t^4+4s^2t^2 & L5a1, L7a1, L7a3, L7a4 \\
12s^4t^4+2s^4t^2+2s^2t^4 & L7n1, \\
12s^4t^4+4s^2t^2 & L6a1, L7a2, L7n2 \\
12s^4t^2+12s^2t^4 +40s^2t^2 & L6a4 \\
16s^2t^2 & L6a5, L6n1, L7a7 \\
\end{array}
\]
In particular we note that most of the links in the table have 16 
$X$-colorings but are distinguished by the enhancement information.
\end{example}

\begin{example}\label{ex:table2}
Let $X$ be the biquandle with operation tables
\[
\begin{array}{r|rrrr}
\utr & 1 & 2 & 3 & 4 \\ \hline
1 & 2 & 2 & 2 & 2 \\
2 & 1 & 1 & 1 & 1 \\
3 & 4 & 4 & 4 & 4 \\
4 & 3 & 3 & 3 & 3
\end{array}
\quad
\begin{array}{r|rrrr}
\otr & 1 & 2 & 3 & 4 \\ \hline
1 & 2 & 2 & 1 & 1 \\
2 & 1 & 1 & 2 & 2 \\
3 & 4 & 4 & 4 & 4 \\
4 & 3 & 3 & 3 & 3
\end{array},
\]
let $\beta$ be the $X$-bracket with values in $R=\mathbb{Z}[q^{\pm 1}]$ given 
by the matrix
\[
\left[\begin{array}{rrrr|rrrr}
1 & 1 & q & q & q & q & 1 & 1 \\
1 & 1 & q & q & q & q & 1 & 1 \\
1 & 1 & q & q & q & q & 1 & 1 \\
1 & 1 & q & q & q & q & 1 & 1 
\end{array}\right].
\]  

The endomorphism ring of $X$ has eight elements; let $S=\mathrm{Hom}(X,X)$. 
Then the trefoil knot $3_1$ has 
\[\Phi_{X,\beta}^{S,\mathrm{deg}^+}(3_1)
=16u^{-q^8 - q^6 - q^4 + 1}v^{12} + 16u^{q^9 - q^5 - q^3 - q}v^4\]
and
\[\Phi_{X,\beta}^{S,2}(3_1)=
16s^{-q^8 - q^6 - q^4 + 1}t^{-q^8 - q^6 - q^4 + 1} + 8s^{q^9 - q^5 - q^3 - q}t^{-q^8 - q^6 - q^4 + 1} + 8s^{q^9 - q^5 - q^3 - q}t^{q^9 - q^5 - q^3 - q}
\]
while the figure eight knot $4_1$ has
\[\Phi_{X,\beta}^{S,\mathrm{deg}^+}(4_1)=16u^{-q^5 -q^{-5}}v^{12} + 16u^{-q^5 -q^{-5}}v^4
\]
and
\[\Phi_{X,\beta}^{S,2}(4_1)=32s^{-q^5 -q^{-5}}t^{-q^5 -q^{-5}}.\]
These two knots have the same biquandle bracket invariant value but are 
distinguished by the quiver enhancement information.
\end{example}

\begin{example}
The table contains the values of $\Phi_{X,\beta}^{S,\mathrm{deg}^+}(L)$ for all 
prime knots with up to eight crossings with respect to the biquandle and
biquandle bracket in Example \ref{ex:table2} and $S=\mathrm{Hom}(X,X)$.
\[
\begin{array}{r|l}
K & \Phi_{X,\beta}^{S,\mathrm{deg}^+}(L) \\ \hline
3_1 & 16u^{-q^8 - q^6 - q^4 + 1}v^{12} + 16u^{q^9 - q^5 - q^3 - q}v^4 \\ 
4_1 & 16u^{-q^5 -q^{-5}}v^{12} + 16u^{-q^5 -q^{-5}}v^4 \\
5_1 & 16u^{-q^{12} - q^{10} - q^8 + 1}v^{12} + 16u^{q^{15} - q^7 - q^5 - q^3}v^4 \\
5_2 & 16u^{-q^{14} - q^{10} - q^8 + q^2}v^{12} + 16u^{q^{13} - q^7 - q^5 - q}v^4\\
6_1 & 16u^{-q^3 + q^{-3} - q^{-7} - q^{-11}}v^{12}
+ 16u^{-q^{-4} - q^{-8} + q^{-12} - q^{-9}}v^4 \\
6_2 & 16u^{-q^2 + 1 - q^{-8} - q^{-9}}v^{12}
+ 16u^{-q^3 - q^{-1} + q^{-9} - q^{-11}}v^4 \\
6_3 & 16u^{q^7 - q^5 - q - q^{-1} - q^{-5} + q^{-7}}v^{12}
+ 16u^{q^7 - q^5 - q - q^{-1} - q^{-5} + q^{-7}}v^4 \\
7_1 & 16u^{-q^{16} - q^{14} - q^{12} + 1}v^{12}
+ 16u^{q^{21} - q^9 - q^7 - q^5}v^4 \\
7_2 & 16u^{q^{-4} - q^{-10} - q^{-16} - q^{-20}}v^{12}
+ 16u^{-q^{-1} - q^{-5} - q^{-11} + q^{-17}}v^4 \\
7_3 & 16u^{-q^{18} - q^{14} - q^{10} - q^8 + q^6 + q^2}v^{12}
+ 16u^{q^{19} + q^{15} - q^{13} - q^{11} - q^7 - q^3}v^4 \\
7_4 & 16u^{-q^{20} + q^{18} - q^{16} - q^{14} - q^{12} - q^{10} + q^8 + q^4}v^{12}
+16u^{q^{17} + q^{13} - q^{11} - q^9 - q^7 - q^5 + q^3 - q}v^4 \\
7_5 & 16u^{q^{-2} - q^{-4} + q^{-6}- 2q^{-14} - q^{-18}}v^{12}
+16u^{-q^{-3} - 2q^{-7} + q^{-15} - q^{-17} + q^{-19}}v^4  \\
7_6 & 16u^{q^4 - q^2 + 1 - q^{-2} - q^{-4} - q^{-8} + q^{-10} - q^{-12}}v^{12}
+16u^{-q^3 + q - q^{-1} - q^{-5} - q^{-7} + q^{-9} - q^{-11} + q^{-13}}v^4 \\
7_7 & 16u^{q^{10} - 2q^8 - q^4 + 1 - q^{-2} + q^{-4} - q^{-6}}v^{12} 
+16u^{-q^9 + q^7 - q^5 + q^3 - q^{-1} - 2q^{-5} + q^{-7}}v^4 \\
\end{array}\]
\[
\begin{array}{r|l}
K & \Phi_{X,\beta}^{S,\mathrm{deg}^+}(L) \\ \hline
8_1 & 16u^{-q + q^{-5} - q^{-13} - q^{-17}}v^{12}
+ 16u^{-q^5 - q + q^{-7} - q^{-13}}v^4 \\
8_2 & 16u^{-q^5 + q^3 + q^{-1} - q^{-5} - q^{-9} - q^{-13}}v^{12}
+16u^{-q - q^{-3} - q^{-7} + q^{-11} + q^{-15} - q^{-17}}v^4 \\
8_3 & 16u^{-q^9 - q^5 + q^3 + q^{-3} - q^{-5} - q^{-9}}v^{12}
+ 16u^{-q^9 - q^5 + q^3 + q^{-3} - q^{-5} - q^{-9}}v^4 \\
8_4 & 16u^{-q^{11} + q^9 - q^7 + q^{-1} - q^{-3} - q^{-7}}v^{12}
+16u^{-q^7 - q^3 + q - q^{-7} + q^{-9} - q^{-11}}v^4 \\
8_5 & 16u^{-q^{13} - 2q^9 + q^3 + q - q^{-1} + q^{-3} - q^{-5}}v^{12} 
+16u^{-q^{17} + q^{15} - q^{13} + q^{11} + q^{9} - 2q^3 - q^{-1}}v^4 \\
8_6 & 16u^{-q^3 + q - q^{-1} + q^{-3} + q^{-9} - 2q^{-11} - q^{-15}}v^{12}
+ 16u^{-q^3 - 2q^{-1} + q^{-3} + q^{-9} - q^{-11} + q^{-13} - q^{-15}}v^4 \\
8_7 & 16u^{q^{11} - q^9 - 2q^5 - q^{-1} + q^{-3} - q^{-5} + q^{-7}}v^{12}
+ 16u^{q^{13} - q^{11} + q^9 - q^7 - 2q - q^{-3} + q^{-5}}v^4 \\
8_8 & 16u^{q^{13} - q^{11} + q^9 - 2q^7 - q^5 - q + q^{-1} - q^{-3} + q^{-5}}v^{12}
+16u^{q^{11} - q^9 + q^7 - q^5 - q - 2q^{-1} + q^{-3} - q^{-5} + q^{-7}}v^4 \\
8_{9} & 16u^{-q^9 + q^7 - q^5 + q^3 - q - q^{-1} + q^{-3} - q^{-5} + q^{-7} - q^{-9}}v^{12}
+ 16u^{-q^9 + q^7 - q^5 + q^3 - q - q^{-1} + q^{-3} - q^{-5} + q^{-7} - q^{-9}}v^4 \\
8_{10} & 16u^{q^{11} - q^9 + q^7 - 2q^5 - q^3 - q - q^{-1} + 2q^{-3} - q^{-5} + q^{-7}}v^{12}
+ 16u^{q^{13} - q^{11} + 2q^9 - q^7 - q^5 - q^3 - 2q + q^{-1} - q^{-3} + q^{-5}}v^4 \\
8_{11} & 16u^{-q^3 + q - q^{-1} + 2q^{-3} - q^{-7} - 2q^{-11} + q^{-13} - q^{-15}}v^{12}
+16u^{-q^3 + q - 2q^{-1} - q^{-5} + 2q^{-9} - q^{-11} + q^{-13} - q^{-15}}v^4  \\
8_{12} & 16u^{-q^9 + q^7 - 2q^5 + q^3 + q^{-3} - 2q^{-5} + q^{-7} - q^{-9}}v^{12} 
+ 16u^{-q^9 + q^7 - 2q^5 + q^3 + q^{-3} - 2q^{-5} + q^{-7} - q^{-9}}v^4 \\
8_{13} & 16u^{q^{13} - 2q^{11} + q^9 - q^7 - 2q + q^{-1} - q^{-3} + q^{-5}}v^{12} 
+16u^{q^{11} - q^9 + q^7 - 2q^5 - q^{-1} + q^{-3} - 2q^{-5} + q^{-7}}v^4  \\
8_{14} & 16u^{-q^3 + 2q - q^{-1} + q^{-3} - q^{-5} - q^{-7}+ q^{-9} - 2q^{-11} + q^{-13} - q^{-15}}v^{12}+16u^{-q^3 + q - 2q^{-1} + q^{-3} - q^{-5} - q^{-7} + q^{-9} -q^{-11} + 2q^{-13} - q^{-15}}v^4  \\
8_{15} & 16u^{-q^{-3} + 2q^{-5} - q^{-7} + 2q^{-9} - q^{-13} - 3q^{-17} + q^{-19} - q^{-21}}v^{12}
+16u^{-q^{-3} + q^{-5} - 3q^{-7} - q^{-11} + 2q^{-15} - q^{-17} + 2q^{-19} - q^{-21}}v^4 
\\
8_{16} & 16u^{q^7 - 2q^5 + 2q^3 - q - 2q^{-5} + q^{-7} - 2q^{-9} + q^{-11}}v^{12}
+16u^{q^5 - 2q^3 + q - 2q^{-1} - q^{-7} + 2q^{-9} - 2q^{-11} + q^{-13}}v^4 \\
8_{17} & 16u^{-q^9 + 2q^7 - 2q^5 + q^3 - q - q^{-1} + q^{-3} - 2q^{-5} + 2q^{-7} - q^{-9}}v^{12}
+ 16u^{-q^9 + 2q^7 - 2q^5 + q^3 - q - q^{-1} + q^{-3} - 2q^{-5} + 2q^{-7} - q^{-9}}v^4 \\
8_{18} & 16u^{-q^9 + 3q^7 - 2.q^5 + q^3 - 2q - 2q^{-1} + q^{-3} - q^{-5} + 3q^{-7} - q^{-9}}v^{12}
+ 16u^{-q^9 + 3q^7 - 2q^5 + q^3 - 2q - 2q^{-1} + q^{-3} - 2q^{-5} + 3q^{-7} - q^{-9}}v^4 \\
8_{19} & 16u^{-q^{19} - q^{17} - q^{15} - q^{13} + q^9 + q^7}v^{12}
+16u^{q^{17} + q^{15} - q^{11} - q^{9} - q^7 - q^5}v^4 \\
8_{20} & 16u^{q^5 -q^{-1} -q^{-3} - q^{-5} - q^{-7} + q^{-9}}v^{12}
+16u^{q^3 - q -q^{-1} -q^{-3} - q^{-5} + q^{-11}}v^4 \\
8_{21} & 16u^{-q^3+q+q^{-3} - q^{-7} - 2q^{-11}}v^{12}
+16u^{-2q^{-1} -q^{-5} + q^{-9} + q^{-13} -q^{-15}}v^4 
\end{array}\]

\end{example}

\begin{example}
Let $X$ be the biquandle with operation tables
\[
\begin{array}{r|rrrr}
\utr & 1 & 2 & 3 & 4 \\ \hline
1 & 2 & 2 & 2 & 2 \\
2 & 1 & 1 & 1 & 1 \\
3 & 4 & 4 & 4 & 4 \\
4 & 3 & 3 & 3 & 3 
\end{array}
\quad
\begin{array}{r|rrrr}
\otr & 1 & 2 & 3 & 4 \\ \hline
1 & 2 & 2 & 1 & 1 \\
2 & 1 & 1 & 2 & 2 \\
3 & 4 & 4 & 4 & 4 \\
4 & 3 & 3 & 3 & 3 
\end{array},
\]
let $\beta$ be the $X$-bracket with values in $R=\mathbb{Z}[q^{\pm 1}]$ 
given by the matrix
\[
\left[\begin{array}{rrrr|rrrr}
1 & 1 & q & q & q & q & 1 & 1 \\
1 & 1 & q & q & q & q & 1 & 1 \\
1 & 1 & q & q & q & q & 1 & 1 \\
1 & 1 & q & q & q & q & 1 & 1 
\end{array}\right]
\]
and let $S=\mathrm{Hom}(X,X)$.
The virtual knots $3.1$ and $3.7$
\[\includegraphics{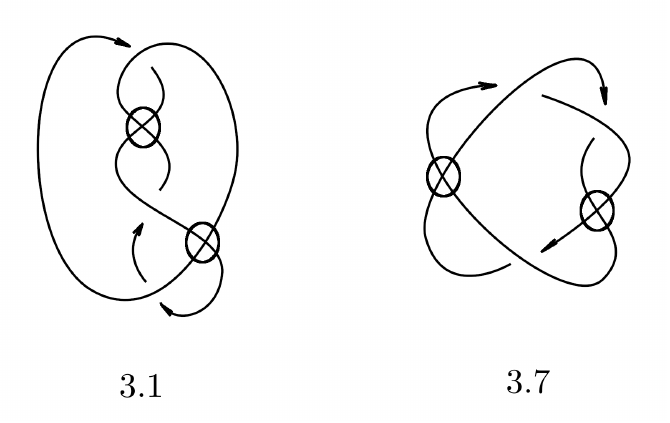}\]
 both have trivial Jones polynomial but are
distinguished by the biquandle bracket quiver invariants with
\[\Phi_{X,\beta}^{S,\mathrm{deg}^+}(3.1)=16u^{-q-q^{-1}}v^4.\]
and
\[\Phi_{X,\beta}^{S,\mathrm{deg}^+}(3.7)=8u^{-q-q^{-1}}v^4+8u^{-q^{-2}-q^{-4}}v^4.\]
In particular, this example shows that biquandle bracket quivers and their
decategorification invariants are not determined by the Jones polynomial.
\end{example}

\section{\large\textbf{Questions}}\label{Q}

We conclude with some questions for future research.

\begin{itemize}
\item One of the more exciting aspects of this project was the realization
that while many biquandle brackets are cohomologous to classical skein 
invariants such as the Jones or HOMFLYPT polynomial, the biquandle bracket
quivers for these can be different. In particular, this means our construction
gives new and different categorifications of these invariants, allowing for 
further new invariants via decategorification.
\item What other new invariants can be obtained from these quivers?
\item How can these quivers be further enhanced?
\item As always, faster method of computing biquandle brackets, particularly
for large biquandles and large finite or infinite rings, are of great interest.
\end{itemize}

\bibliography{pf-sn}{}

\begin{thebibliography}{1}

\bibitem{KA}
D.~Bar-Natan.
\newblock The knot atlas \textup{http://katlas.org/wiki/Main\_Page}.

\bibitem{CES}
J.~S. Carter, M.~Elhamdadi, and M.~Saito.
\newblock Homology theory for the set-theoretic {Y}ang-{B}axter equation and
  knot invariants from generalizations of quandles.
\newblock {\em Fund. Math.}, 184:31--54, 2004.

\bibitem{CCN}
J.~Ceniceros, A.~Christiana, and S.~Nelson.
\newblock Psyquandle coloring quivers.
\newblock {\em arXiv:2107.05668}, 2021.

\bibitem{CN}
K.~Cho and S.~Nelson.
\newblock Quandle coloring quivers.
\newblock {\em Journal of Knot Theory and Its Ramifications}, 28(01):1950001,
  2019.

\bibitem{EN}
M.~Elhamdadi and S.~Nelson.
\newblock {\em Quandles---an introduction to the algebra of knots}, volume~74
  of {\em Student Mathematical Library}.
\newblock American Mathematical Society, Providence, RI, 2015.

\bibitem{HVW}
W.~Hoffer, A.~Vengal, and V.~Winstein.
\newblock The structure of biquandle brackets.
\newblock {\em J. Knot Theory Ramifications}, 29(6):2050042, 13, 2020.

\bibitem{KR}
L.~H. Kauffman and D.~Radford.
\newblock Bi-oriented quantum algebras, and a generalized {A}lexander
  polynomial for virtual links.
\newblock In {\em Diagrammatic morphisms and applications ({S}an {F}rancisco,
  {CA}, 2000)}, volume 318 of {\em Contemp. Math.}, pages 113--140. Amer. Math.
  Soc., Providence, RI, 2003.

\bibitem{NOR}
S.~Nelson, M.~E. Orrison, and V.~Rivera.
\newblock Quantum enhancements and biquandle brackets.
\newblock {\em J. Knot Theory Ramifications}, 26(5):1750034, 24, 2017.

\end{thebibliography}
\bibliographystyle{abbrv}

\bigskip

\noindent
\textsc{Department of Mathematical Sciences \\
Claremont McKenna College \\
850 Columbia Ave. \\
Claremont, CA 91711}

\end{document}